\documentclass[reqno, 11pt, a4paper]{article}
\usepackage{amssymb,dsfont,graphicx}
\usepackage{amsthm,amsmath}
\usepackage[usenames,dvipsnames]{color}
\usepackage{fullpage}

\ifx\pdfoutput\undefined%
  \newcommand\phantomsection\relax
  \newcommand{\url}[1]{#1}
  \newcommand{\href}[2]{#2}
\else
\usepackage[pdftex,colorlinks]{hyperref}
\hypersetup{%
pdftitle={On queues with service and interarrival times depending on waiting times},
pdfauthor={O.J. Boxma, M. Vlasiou},
pdfkeywords={},
bookmarksnumbered,
pdfstartview={FitH},
citecolor=Plum,
linkcolor=RedOrange,
urlcolor=RawSienna,
}%
\fi

\theoremstyle{plain}              
\newtheorem{theorem}{Theorem}
\newtheorem{lemma}{Lemma}
\theoremstyle{remark}
\newtheorem{remark}{Remark}

\numberwithin{equation}{section}    

\newcommand{\m}[1]{\mathcal{#1}}
\newcommand{\e}{\mathbb{E}}
\newcommand{\p}{\mathbb{P}}

\newcommand{\Dfb}[1][]{\mbox{$F_B^{#1}$}}
\newcommand{\Dfa}[1][]{\mbox{$F_A^{#1}$}}
\newcommand{\Dfw}[1][]{\mbox{$F_W^{#1}$}}
\newcommand{\dfw}[1][]{\mbox{$f_W^{#1}$}}
\newcommand{\Dfx}[1][]{\mbox{$F_X^{#1}$}}

\newcommand{\lta}[1][]{\mbox{$\alpha^{#1}$}}
\newcommand{\ltw}[1][]{\mbox{$\omega^{#1}$}}

\begin{document}
\title{On queues with service and interarrival times depending on waiting times}
\author{O.J.\ Boxma$^{\star,\,\star\star}$, M.\ Vlasiou$^{\star\star\star,\,}$\footnote{This research has been carried out when this author was affiliated with EURANDOM, The Netherlands.}}
\date{November 26, 2006}
\maketitle
\begin{center}
$^\star$ EURANDOM,\\
P.O.\ Box 513, 5600 MB Eindhoven, The Netherlands.
\end{center}

\begin{center}
$^{\star\star}$ Eindhoven University of Technology,
\\Department of Mathematics \& Computer Science,
\\P.O.\ Box 513, 5600 MB Eindhoven, The Netherlands.
\end{center}

\begin{center}
$^{\star\star\star}$ Georgia Institute of Technology,
\\H. Milton Stewart School of Industrial \& Systems Engineering,
\\765 Ferst Drive, Atlanta GA 30332-0205, USA.
\\\vspace{0.3cm}\href{mailto:boxma@win.tue.nl}{boxma@win.tue.nl}, \href{mailto:vlasiou@gatech.edu}{vlasiou@gatech.edu}
\end{center}

\begin{abstract}
We consider an extension of the standard $\mathrm{G/G/1}$ queue, described by the equation $W\stackrel{\m{D}}{=}\max\{0, B-A+YW\}$, where $\p[Y=1]=p$ and $\p[Y=-1]=1-p$. For $p=1$ this model reduces to the classical Lindley equation for the waiting time in the $\mathrm{G/G/1}$ queue, whereas for $p=0$ it describes the waiting time of the server in an alternating service model. For all other values of $p$ this model describes a FCFS queue in which the service times and interarrival times depend linearly and randomly on the waiting times. We derive the distribution of $W$ when $A$ is generally distributed and $B$ follows a phase-type distribution, and when $A$ is exponentially distributed and $B$ deterministic.
\end{abstract}

\section{Introduction}\label{s:intro}
One of the most fundamental relations in queuing and random walk theory is Lindley's recursion \cite{lindley52}:
\begin{equation}\label{eq:classic}
 W_{n+1}=\max\{0, B_n-A_n+W_n\}.
\end{equation}
In queuing theory, it represents a relation between the waiting times of the $n$-th and $(n+1)$-st customer in a single server queue, $A_n$ indicating the interarrival time between the $n$-th and $(n+1)$-st customer and $B_n$ denoting the service time of the $n$-th customer. In the applied probability literature there has been a considerable amount of interest in generalisations of Lindley's recursion, namely the class of Markov chains described by the recursion $W_{n+1}=g(W_n,X_n)$. For earlier work on such stochastic recursions see for example Brandt {\em et al.}~\cite{brandt-SSM} and Borovkov and Foss~\cite{borovkov92}. Many structural properties of this recursion have been derived. For example, Asmussen and Sigman~\cite{asmussen96a} develop a duality theory, relating the steady-state distribution to a ruin probability associated with a risk process. More references in this domain can be found for example in Asmussen and Schock Petersen~\cite{asmussen89} and Seal~\cite{seal72}. An important assumption which is often made in these studies is that the function $g(w,x)$ is non-decreasing in its main argument $w$. For example, in \cite{asmussen96a} this assumption is crucial for their duality theory to hold.

In this paper we consider a generalisation of Lindley's recursion, for which the monotonicity assumption does not hold. In particular, we study the Lindley-type recursion
\begin{equation}\label{eq:recursion}
 W_{n+1}=\max\{0, B_n-A_n+Y_nW_n\},
\end{equation}
where for every $n$, the random variable $Y_n$ is equal to plus or minus one according to the probabilities $\p[Y_n=1]=p$ and $\p[Y_n=-1]=1-p$, $0 \leqslant p\leqslant 1$. The sequences $\{A_n\}$ and $\{B_n\}$ are assumed to be independent sequences of i.i.d.\ non-negative random variables. Our main goal is to derive the steady-state distribution of $\{W_n,~n=1,2,\dots\}$, when it exists.

Equation \eqref{eq:recursion} reduces to the classical Lindley recursion \cite{lindley52} when $\p[Y_n=1]=1$ for every $n$. Furthermore, if $\p[Y_n=-1]=1$, then \eqref{eq:recursion} describes the waiting time of the server in an alternating service model with two service points. For a description of the model for this special case and related results see Park {\em et al.}~\cite{park03} and Vlasiou {\em et al.}~\cite{vlasiou05a, vlasiou05, vlasiou05b, vlasiou04}.

Studying a recursion that contains both Lindley's classical recursion and the recursion in \cite{park03, vlasiou05a, vlasiou05, vlasiou05b, vlasiou04} as special cases seems of interest in its own right. Additional motivation for studying the recursion is supplied by the fact that, for $0<p<1$, the resulting model can be interpreted as a special case of a queuing model in which service and interarrival times depend on waiting times. We shall now discuss the latter model.

Consider an extension of the standard $\mathrm{G/G/1}$ queue in which the service times and the interarrival times depend linearly and randomly on the waiting times. Namely, the model is specified by a stationary and ergodic sequence of four-tuples of nonnegative random variables $\{(A_n, B_n, \widehat{A}_n, \widehat{B}_n)\}$, $n\geqslant 0$. The sequence $\{W_n\}$ is defined recursively by
$$
 W_{n+1}=\max\{0, \overline{B}_n-\overline{A}_n+W_n\},
$$
where
\begin{align*}
  \overline{A}_n&=A_n+\widehat{A}_nW_n,\\
  \overline{B}_n&=B_n+\widehat{B}_nW_n.
\end{align*}

We interpret $W_n$ as the waiting time and $\overline{B}_n$ as the service time of customer $n$. Furthermore, we take $\overline{A}_n$ to be the interarrival time between customers $n$ and $n+1$. We call $B_n$ the \textit{nominal service time} of customer $n$ and $A_n$ the \textit{nominal interarrival time} between customers $n$ and $n+1$, because these would be the actual times if the additional shift were omitted, that is, if $\p[\widehat{A}_n=\widehat{B}_n=0]=1$.

Evidently, the waiting times satisfy the generalised Lindley recursion \eqref{eq:recursion}, where we have written $Y_n=1+\widehat{B}_n-\widehat{A}_n$. This model -- for generally distributed random variables $Y_n$ -- has been introduced in Whitt~\cite{whitt90}, where the focus is on conditions for the process to converge to a proper steady-state limit, and on approximations for this limit. There are very few exact results known for queuing models in which interarrival and/or service times depend on waiting times; we refer to Whitt \cite{whitt90} for some references.

Whitt \cite{whitt90} builds upon previous results by Vervaat~\cite{vervaat79} and Brandt~\cite{brandt86} for the unrestricted recursion $W_{n+1}=Y_nW_n+X_n$, where $X_n=B_n-A_n$.
There has been considerable previous work on this unrestricted recursion, due to its close connection to the problem of the ruin of an insurer who is exposed to a stochastic economic environment. Such an environment has two kinds of risk, which were called by Norberg~\cite{norberg99} insurance risk and financial risk. Indicatively, we mention the work by Tang and Tsitsiashvili~\cite{tang03}, and by Kalashnikov and Norberg~\cite{kalashnikov02a}. In the more general framework, $W_n$ may represent an inventory in time period $n$ (e.g.\ cash), $Y_n$ may represent a multiplicative, possibly random, decay or growth factor between times $n$ and $n+1$ (e.g.\ interest rate) and $B_n-A_n$ may represent a quantity that is added or subtracted between times $n$ and $n+1$ (e.g.\ deposit minus withdrawal). Obviously, the positive-part operator is appropriate for many applications~\cite{whitt90}.

This paper presents an exact analysis of the steady-state distribution of $\{W_n, ~ n=1,2,\dots\}$ as given by
\eqref{eq:recursion} with $\p[Y_n=1]=p$ and $\p[Y_n=-1]=1-p$. For $0<p<1$, this amounts to analysing
the above-described $\mathrm{G/G/1}$ extension where $\widehat{A}_n=\widehat{B}_n$ with probability $p$, and
$\widehat{A}_n=2+\widehat{B}_n$ with probability $1-p$. This problem, and state-dependent queuing processes in general, is connected to LaPalice queuing models, introduced by Jacquet~\cite{jacquet92}, where customers are scheduled in such a way that the period between two consecutively scheduled customers is greater than or equal to the service time of the first customer.

This paper is organised in the following way. In Section~\ref{s:model} we comment on the stability of the process $\{W_n\}$, as it is defined by Recursion \eqref{eq:recursion}. In the remainder of the paper it is assumed that the steady-state distribution of $\{W_n\}$ exists. Section~\ref{s:G/PH} is devoted to the determination of the distribution of $W$ when $A$ is generally distributed and $B$ has a phase-type distribution. In Section~\ref{s:M/D} we determine the distribution of $W$ when $A$ is exponentially distributed and $B$ is deterministic. At the end of each section we compare the results that we derive to the already known results for Lindley's recursion (i.e.\ for $p=1$) and to the equivalent results for the Lindley-type recursion arising for $p=0$.

At the end of this introduction we mention a few notational conventions. For a random variable $X$ we denote its distribution by $F_X$ and its density by $f_X$. Furthermore, we shall denote by $f^{(i)}$ the $i$-th derivative of the function $f$. The Laplace-Stieltjes transforms (LST) of $A$ and $W$ are respectively denoted by $\lta$ and $\ltw$. To keep expressions simple, we also use the function $\phi$ defined as $\phi(s)=\omega(s) \, \lta(s)$.

\section{Stability}\label{s:model}
The following result on the convergence of the process $\{W_n\}$ to a proper limit $W$ is shown in Whitt~\cite{whitt90}. It is included here only for completeness.

From Recursion~\eqref{eq:recursion}, it is obvious that if we replace $Y_n$ by $\max\{0, Y_n\}$ and $B_n-A_n$ by $\max\{0, B_n-A_n\}$, then the resulting waiting times will be at least as large as the ones given by \eqref{eq:recursion}. Moreover, when we make this change, the positive-part operator is not necessary anymore.
\begin{lemma}[Whitt~\protect{\cite[Lemma 1]{whitt90}}]
If\/ $W_n$ satisfies \eqref{eq:recursion}, then with probability 1, $W_n\leqslant Z_n$ for all $n$, where
\begin{equation}\label{eq:unrestricted}
Z_{n+1}=\max\{0, Y_n\}Z_{n}+\max\{0, B_n-A_n\},\quad n\geqslant 0,
\end{equation}
and $Z_0=W_0\geqslant 0$.
\end{lemma}
So if $W_n$ satisfies \eqref{eq:recursion}, $Z_n$ satisfies \eqref{eq:unrestricted}, and $Z_n$ converges to the proper limit $Z$, then $\{W_n\}$ is tight and $\p[W>x]\leqslant\p[Z>x]$ for all $x$, where $W$ is the limit in distribution of any convergent subsequence of $\{W_n\}$. This observation, combined with Theorem~1 of Brandt~\cite{brandt86}, which implies that $Z_n$ satisfying \eqref{eq:unrestricted} converges to a proper limit if $\p[\max\{0, Y_n\}=0]=\p[Y_n\leqslant 0]>0$, leads to the following theorem.
\begin{theorem}[Whitt~\protect{\cite[Theorem 1]{whitt90}}]\label{th:stability}
The series $\{W_n\}$ is tight for all $\rho=\e[B_0]/\e[A_0]$ and $W_0$. If, in addition, $0\leqslant p<1$ and $\{(Y_n, B_n-A_n)\}$ is a sequence of independent vectors with $$\p[Y_0\leqslant0, B_0-A_0\leqslant 0]>0,$$ then the events $\{W_n=0\}$ are regeneration points with finite mean time and $\{W_n\}$ converges in distribution to a proper limit $W$ as $n\to\infty$ for all $\rho$ and $W_0$.
\end{theorem}
Naturally, for $p=1$, i.e.\ for the classical Lindley recursion, we need the additional condition that $\rho <1$.

Therefore, assume that the sequences $\widehat{B}_n-\widehat{A}_n$ and $B_n-A_n$ are independent stationary sequences, that are also independent of one another, and that for all $n$, $A_n$ and $B_n$ are non-negative. Then the conditions of Theorem~\ref{th:stability} hold, so there exists a proper limit $W$, and for the system in steady-state we write
\begin{equation}\label{eq:equation}
W\stackrel{\m{D}}{=}\max\{0, B-A+YW\},
\end{equation}
where ``$\stackrel{\m{D}}{=}$'' denotes equality in distribution, where $A$, $B$ are generic random variables distributed like $A_n$, $B_n$, and where $\p[Y=1]=p$ and $\p[Y=-1]=1-p$.

\begin{remark}
For $x\geqslant 0$ Equation~\eqref{eq:equation} yields that
$$
\Dfw(x)=\p[W\leqslant x]=p\, \p[X+W\leqslant x]+(1-p) \p[X-W\leqslant x],
$$
where $X=B-A$ (note that $\p[X<0]>0$). Assuming that the distribution $\Dfx$ of the random variable $X$ is continuous, the last term is equal to $1-\p[X-W\geqslant x]$, which gives us that
$$
\Dfw(x)=p \int^x_{-\infty}\Dfw(x-y)\,{\rm d}\Dfx(y)+(1-p)\biggl(1-\int_x^\infty\Dfw(y-x)\,{\rm d}\Dfx(y)\biggr).
$$
This means that the limiting distribution of $W$, provided that $\Dfx$ is continuous, satisfies the functional equation
\begin{equation}\label{eq:funtional eq}
F(x)=p \int^x_{-\infty} F(x-y)\,{\rm d}\Dfx(y)+(1-p)\biggl(1-\int_x^\infty F(y-x)\,{\rm d}\Dfx(y)\biggr).
\end{equation}
Therefore, there exists at least one function that is a solution to \eqref{eq:funtional eq}. It can be shown that in fact there exists a unique measurable bounded function $F: [0,\infty) \to \mathds{R}$ that satisfies this functional equation.

To show this, consider the space $\mathcal{L}^\infty([0,\infty))$, i.e.\ the space of measurable and bounded functions on the real line with the norm
$$
\|F\|= \sup_{t \geqslant 0} |F(t)|.
$$
In this space define the mapping
$$
(\mathcal{T}F)(x)=p \int^x_{-\infty} F(x-y)\,{\rm d}\Dfx(y)+(1-p)\biggl(1-\int_x^\infty F(y-x)\,{\rm d}\Dfx(y)\biggr).
$$
Note that $\mathcal{T}F:\mathcal{L}^\infty\bigl([0,\infty)\bigr) \rightarrow \mathcal{L}^\infty\bigl([0,\infty)\bigr)$, i.e., $\mathcal{T}F$ is measurable and bounded. For two arbitrary functions $F_1$ and $F_2$ in this space we have
\begin{align*}
\|&(\mathcal{T}F_1)-(\mathcal{T}F_2)\|=\sup_{x \geqslant 0} \left|(\mathcal{T}F_1)(x)-(\mathcal{T}F_2)(x)\right|\\
&=\sup_{x \geqslant 0}\ \left|p\int^x_{-\infty}\left[F_1(x-y)-F_2(x-y)\right]\,{\rm d}\Dfx(y)+(1-p)\int_x^\infty\left[F_2(y-x)-F_1(y-x)\right]\,{\rm d}\Dfx(y)\right|\\
&\leqslant \sup_{x \geqslant 0}\left(p \int^x_{-\infty}|F_1(x-y)-F_2(x-y)|\,{\rm d}\Dfx(y)+(1-p)\int_x^\infty |F_2(y-x)-F_1(y-x)|\,{\rm d}\Dfx(y)\right)\\
&\leqslant \sup_{x \geqslant 0}\left(p \int^x_{-\infty}\sup_{t \geqslant 0}|F_1(t)-F_2(t)|\,{\rm d}\Dfx(y)+(1-p)\int_x^\infty \sup_{t \geqslant 0}|F_2(t)-F_1(t)|\,{\rm d}\Dfx(y)\right)\\
&=\|F_1-F_2\|\ \sup_{x \geqslant 0}\bigl[p\,\Dfx(x)+(1-p)\bigr(1-\Dfx(x)\bigr)\bigr].
\end{align*}
Note that the supremum appearing above is less than or equal to $\max\{p,1-p\}$ for $p\in(0,1)$ and equal to $1-\Dfx(0)$ for $p=0$. Therefore, since for $p\neq 1$ it holds that $\sup_{x \geqslant 0}\bigl[p\,\Dfx(x)+(1-p)\bigr(1-\Dfx(x)\bigr)\bigr]<1$, for these values of the parameter $p$ we have a contraction mapping. Furthermore, we know that $\mathcal{L}^\infty([0,\infty))$ is a Banach space, therefore by the Fixed Point Theorem we have that \eqref{eq:funtional eq} has a unique solution; for a related discussion, see also \cite{vlasiou05a}.
\end{remark}

In the next two sections we determine the distribution of $W$ for two cases in which $\{A_n\}$ and $\{B_n\}$ are independent and i.i.d.\ sequences of non-negative random variables.

\section{The GI/PH case}\label{s:G/PH}
In this section we assume that $A$ is generally distributed, while $B$ follows a particular phase-type distribution. Specifically, we assume that with probability $\kappa_n$ the nominal service time $B$ follows  an Erlang distribution with parameter $\mu$ and $n$ phases, i.e.,
\begin{equation}\label{eq:distr B}
\Dfb(x)=\sum_{n=1}^N \kappa_n\biggl(1-\mathrm{e}^{-\mu x}\sum_{j=0}^{n-1}\frac{(\mu x)^j}{j!}\biggr)=\sum_{n=1}^N \kappa_n\sum_{j=n}^{\infty} \mathrm{e}^{-\mu x}\frac{(\mu x)^j}{j!},\qquad x \geqslant 0,
\end{equation}
with LST $\sum_{n=1}^N \kappa_n \bigl(\mu/(\mu+s)\bigr)^n$. These distributions, i.e.\ mixtures of Erlang distributions, are special cases of Coxian or phase-type distributions. It is sufficient to consider only this class, since it may be used to approximate any given continuous distribution on $[0,\infty)$ arbitrarily close; see Schassberger~\cite{schassberger-W}. Following the proof in \cite{vlasiou05b}, one can show that for such an approximation of $\Dfb$, the error in the resulting waiting time approximation can be bounded.

We are interested in the distribution of $W$. In order to derive the distribution of $W$, we shall first derive the LST of $\Dfw$. We follow a method based on Wiener-Hopf decomposition. A straightforward calculation yields for values of $s$ such that $\mathrm{Re}(s)=0$:
\begin{align}
\nonumber \ltw(s)&=\e[\mathrm{e}^{-sW}]=p\,\e[\mathrm{e}^{-s\max\{0, B-A+W\}}]+(1-p)\e[\mathrm{e}^{-s\max\{0, B-A-W\}}]\\
\nonumber        &=p\,\p[W+B\leqslant A]+p\,\e[\mathrm{e}^{-s(B-A+W)}]-p\,\e[\mathrm{e}^{-s(B-A+W)};W+B\leqslant A]+\\
\label{eq:ltw1}  &\quad +(1-p) \p[B\leqslant W+A]+(1-p)\e[\mathrm{e}^{-s(B-A-W)};B\geqslant W+A];
\end{align}
here $A$, $B$ and $W$ are independent random variables. The Lindley-type equation $W\stackrel{\m{D}}{=}\max\{0, B-A-W\}$ for $A$ generally distributed and $B$ phase-type has already been analysed in Vlasiou and Adan~\cite{vlasiou05}, and the LST of the corresponding $W$ is given there. From Equation~(3.8) of \cite{vlasiou05} we can readily copy an expression for the last two terms appearing in \eqref{eq:ltw1}, so $\ltw$ can now be written as
\begin{align*}
\ltw(s)&=p\,\p[W+B\leqslant A]+p\,\lta(-s)\ltw(s)\sum_{n=1}^N \kappa_n\biggl(\frac{\mu}{\mu+s}\biggr)^n -p\,\e[\mathrm{e}^{-s(B-A+W)};W+B \leqslant A]+\\
        &\quad +(1-p)\Biggl[1-\sum_{n=1}^N \sum_{i=0}^{n-1} \kappa_n \frac{(-\mu)^i}{i!} \phi^{(i)}(\mu)\biggl(1-\biggl( \frac{\mu}{\mu+s}\biggr)^{n-i} \biggr)\Biggr].
\end{align*}
So for $\mathrm{Re}(s)=0$ we have that
\begin{multline}\label{eq:transform1}
\ltw(s)\Biggl[1-p\,\lta(-s)\sum_{n=1}^N \kappa_n\biggl(\frac{\mu}{\mu+s}\biggr)^n  \Biggr]=p\,\p[W+B\leqslant A]-p\,\e[\mathrm{e}^{-s(B-A+W)};W+B\leqslant A]+\\
        +(1-p)\Biggl[1-\sum_{n=1}^N \sum_{i=0}^{n-1} \kappa_n \frac{(-\mu)^i}{i!}\, \phi^{(i)}(\mu) \biggl(1 - \biggl(\frac{\mu}{\mu+s}\biggr)^{n-i} \biggr)\Biggr].
\end{multline}
Cohen~\cite[p.\ 322--323]{cohen-SSQ} shows by applying Rouch\'{e}'s theorem that the function
$$
1-p\,\lta(-s)\sum_{n=1}^N \kappa_n\biggl(\frac{\mu}{\mu+s}\biggr)^n \equiv \frac{1}{(\mu+s)^N}\biggl[(\mu+s)^N-p\,\lta(-s)\sum_{n=1}^N \kappa_n \mu^n(\mu+s)^{N-n}\biggr]
$$
has exactly $N$ zeros $\xi_i(p)$ in the left-half plane if $0<p<1$ (it is assumed that $\lta(\mu) \neq 0$, which is not an essential restriction) or if $p=1$ and $\e[B]<\e[A]$. Naturally, this statement is not valid if $p=0$; therefore, this case needs to be excluded from this point on. So we rewrite \eqref{eq:transform1} as follows
\begin{multline}\label{eq:transform2}
\ltw(s)\prod_{i=1}^{N}\bigl(s-\xi_i(p)\bigr)=\frac{\prod_{i=1}^{N}\bigl(s-\xi_i(p)\bigr)}{(\mu+s)^N-p\,\lta(-s)\sum_{n=1}^N \kappa_n \mu^n(\mu+s)^{N-n}} \times\\
\times\Biggl[p\,(\mu+s)^N\p[W+B\leqslant A]-p\,(\mu+s)^N \e[\mathrm{e}^{-s(B-A+W)};W+B\leqslant A]+\\
+(1-p)\biggl[(\mu+s)^N-\sum_{n=1}^N \sum_{i=0}^{n-1} \kappa_n \frac{(-\mu)^i}{i!} \phi^{(i)}(\mu) \bigl((\mu+s)^N- \mu^{n-i}(\mu+s)^{N-n+i} \bigr)\biggr]\Biggr].
\end{multline}
The left-hand side of \eqref{eq:transform2} is analytic for $\mathrm{Re}(s)>0$ and continuous for $\mathrm{Re}(s)\geqslant0$, and the right-hand side of \eqref{eq:transform2} is analytic for $\mathrm{Re}(s)<0$ and continuous for $\mathrm{Re}(s)\leqslant0$. So from Liouville's theorem \cite{titchmarsh-tf} we have that both sides of \eqref{eq:transform2} are the same $N$-th degree polynomial, say, $\sum_{i=0}^N q_i s^i$. Hence,
\begin{equation}\label{eq:transform3}
  \ltw(s)=\frac{\sum_{i=0}^N q_i s^i}{\prod_{i=1}^{N}\bigl(s-\xi_i(p)\bigr)}.
\end{equation}
In the expression above, the constants $q_i$ are not determined so far, while the roots $\xi_i(p)$ are known. In order to obtain the transform, observe that $\ltw$ is a fraction of two polynomials of degree $N$. So, ignoring the special case of multiple zeros $\xi_i(p)$, partial fraction decomposition yields that \eqref{eq:transform3} can be rewritten as
\begin{equation}\label{eq:transform4}
  \ltw(s)=c_0+ \sum_{i=1}^N \frac{c_i}{\bigl(s-\xi_i(p)\bigr)},
\end{equation}
which implies that the waiting time distribution has a mass at the origin that is given by
$$
\p[W=0]=\lim_{s\to\infty}\e[\mathrm{e}^{-sW}]=c_0
$$
and has a density that is given by
$$
\dfw(x)=\sum_{i=1}^Nc_i \mathrm{e}^{\xi_i(p)x}.
$$
All that remains is to determine the $N+1$ constants $c_i$. To do so, we work as follows.

We shall substitute \eqref{eq:transform4} in the left-hand side of \eqref{eq:transform2}, and express the terms $\p[W+B\leqslant A]$ and $\e[\mathrm{e}^{-s(B-A+W)};W+B\leqslant A]$ that appear at the right-hand side of \eqref{eq:transform2} in terms of the constants $c_i$. Note that the terms $\phi^{(i)}(\mu)$ that appear at the right-hand side of \eqref{eq:transform2} can also be expressed in terms of the constants $c_i$. Thus we obtain a new equation that we shall differentiate a total of $N$ times. We shall evaluate each of these derivatives for $s=0$ and thus we obtain a linear system of $N$ equations for the constants $c_i$, $i=0,\ldots,N$. The last equation that is necessary to uniquely determine the constants $c_i$ is the normalisation equation
\begin{equation}\label{norma}
c_0+\int_0^\infty \dfw(x)\,\mathrm{d}x=1.
\end{equation}

To begin with, note that
\begin{equation}\label{eq:w+b}
\p[W+B\leqslant A]=\p[W=0]\p[B\leqslant A]+\int_{0}^\infty\p[B\leqslant A-x] \sum_{i=1}^N c_i \mathrm{e}^{\xi_i(p)x} \,\mathrm{d}x,
\end{equation}
with
\begin{equation}\label{eq:b<a}
\p[B\leqslant A]= \int_{0}^\infty \sum_{n=1}^N \kappa_n \biggl(\mathrm{e}^{-\mu x}\sum_{j=n}^{\infty}\frac{(\mu x)^j}{j!} \biggr)\,\mathrm{d}\Dfa(x) = \sum_{n=1}^N \sum_{i=n}^\infty \kappa_n \frac{(-\mu)^i}{i!} \lta[(i)](\mu),
\end{equation}
and
\begin{align}
\nonumber \int_{0}^\infty\p[B\leqslant A-x]&\sum_{i=1}^N c_i \mathrm{e}^{\xi_i(p)x} \,\mathrm{d}x=\int_0^\infty\int_0^\infty\p[B\leqslant y-x]\sum_{i=1}^N c_i \mathrm{e}^{\xi_i(p)x} \,\mathrm{d}x \mathrm{d}\Dfa(y)\\
\nonumber    &=\int_0^\infty\int_0^y \mathrm{e}^{-\mu(y-x)}\sum_{n=1}^N \sum_{j=n}^\infty\kappa_n\frac{\bigl(\mu(y-x)\bigr)^j}{j!}\sum_{i=1}^N c_i \mathrm{e}^{\xi_i(p)x}\,\mathrm{d}x \mathrm{d}\Dfa(y)\\
\label{eq:int B<A-x}&=\sum_{n=1}^N\sum_{j=n}^\infty\sum_{i=1}^N\sum_{k=j+1}^\infty\kappa_n c_i\frac{\mu^j\bigl(\mu+\xi_i(p)\bigr)^ {k-j-1}}{k!(-1)^k}\lta[(k)](\mu).
\end{align}
Likewise, we have that
\begin{multline}\label{eq:mt}
\e[\mathrm{e}^{-s(B-A+W)};W+B\leqslant A]=\p[W=0]\e[\mathrm{e}^{-s(B-A)};B\leqslant A]+\\
+\int_{0}^\infty\e[\mathrm{e}^{-s(B-A+x)};x+B\leqslant A]\sum_{i=1}^N c_i \mathrm{e}^{\xi_i(p)x} \,\mathrm{d}x,
\end{multline}
with
\begin{align}\label{eq:mt xwris w}
\nonumber \e[\mathrm{e}^{-s(B-A)};B\leqslant A]&=\int_0^\infty\int_0^x \mathrm{e}^{-s(y-x)} \sum_{n=1}^N \kappa_n \mu \mathrm{e}^{-\mu y} \frac{(\mu y)^{n-1}}{(n-1)!}\,\mathrm{d}y \mathrm{d}\Dfa(x) \\
\nonumber &=\int_0^\infty \mathrm{e}^{xs}\sum_{n=1}^N\kappa_n\biggl(\frac{\mu}{\mu+s}\biggr)^n\sum_{i=n}^\infty \mathrm{e}^{-x(\mu+s)}\frac{x^i(\mu+s)^i}{i!} \,\mathrm{d}\Dfa(x) \\
          &= \sum_{n=1}^N \sum_{i=n}^\infty \kappa_n\, {\mu}^n \frac{(-1)^i}{i!}(\mu+s)^{i-n}\lta[(i)](\mu),
\end{align}
and
\begin{align}\label{eq:int mt}
\nonumber \int_{0}^\infty&\e[\mathrm{e}^{-s(B-A+x)};x+B\leqslant A]\sum_{i=1}^N c_i \mathrm{e}^{\xi_i(p)x} \,\mathrm{d}x\\
\nonumber &=\int_0^\infty\int_0^y\int_0^{y-x} \mathrm{e}^{-s(z-y+x)} \sum_{n=1}^N \kappa_n \mu \mathrm{e}^{-\mu z} \frac{(\mu z)^{n-1}}{(n-1)!} \sum_{i=1}^N c_i \mathrm{e}^{\xi_i(p)x}\,\mathrm{d}z \mathrm{d}x \mathrm{d}\Dfa(y)\\
\nonumber &=\int_0^\infty\int_0^y \sum_{n=1}^N \kappa_n \biggl(\frac{\mu}{\mu+s}\biggr)^n \mathrm{e}^{-s(x-y)} \sum_{j=n}^\infty \mathrm{e}^{-(\mu+s)(y-x)} \frac{(\mu+s)^j (y-x)^j}{j!} \sum_{i=1}^N  c_i \mathrm{e}^{\xi_i(p)x}\,\mathrm{d}x \mathrm{d}\Dfa(y)\\
&= \sum_{n=1}^N \sum_{j=n}^\infty\sum_{i=1}^N\sum_{k=j+1}^\infty \kappa_n c_i  \biggl(\frac{\mu}{\mu+s}\biggr)^n \frac{ (\mu+s)^j \bigl(\mu+\xi_i(p)\bigr)^{k-j-1}}{k!(-1)^k}\,\lta[(k)](\mu).
\end{align}

So, using \eqref{eq:b<a} and \eqref{eq:int B<A-x}, substitute \eqref{eq:w+b} in the right-hand side of \eqref{eq:transform2}, and similarly for \eqref{eq:mt}. Furthermore, as mentioned before, substitute \eqref{eq:transform4} into the left-hand side of \eqref{eq:transform2} to obtain an expression, where both sides can be reduced to an $N$-th degree polynomial in $s$. By evaluating this polynomial and all its derivatives for $s=0$ we obtain $N$ equations binding the constants $c_i$. These equations, and the normalisation equation \eqref{norma}, form a linear system for the constants $c_i$, $i=0,\ldots,N$, that uniquely determines them (see also Remark~\ref{rem:complex roots} below). For example, the first equation, evaluated at $s=0$, yields that
$$
c_0-\sum_{i=1}^N\frac{c_i}{\xi_i(p)}=\frac{1-p}{1-p\,\lta(0)}=1,
$$
since $\lta(0)=1$. We summarise the above in the following theorem.
\begin{theorem}\label{th:density}
Consider the recursion given by \eqref{eq:recursion}, and assume that $0<p<1$. Let \eqref{eq:distr B} be the distribution of the random variable $B$. Then the limiting distribution of the waiting time has mass $c_0$ at the origin and a density on $[0, \infty)$ that is given by
\begin{equation*}
\dfw(x)=\sum_{i=1}^N c_i \mathrm{e}^{\xi_i(p)x}.
\end{equation*}
In the above equation, the constants $\xi_i(p)$, with $\mathrm{Re}(\xi_i(p))<0$, are the $N$ roots of
$$
(\mu+s)^N-p\,\lta(-s)\sum_{n=1}^N \kappa_n \mu^n(\mu+s)^{N-n}=0,
$$
and the $N+1$ constants $c_i$ are the unique solution to the linear system described above.
\end{theorem}
\begin{remark}\label{rem:complex roots}
Although the roots $\xi_i(p)$ and coefficients $c_i$ may be complex, the density and the mass $c_0$ at zero will be positive. This follows from the fact that there is a unique equilibrium distribution and thus a unique solution to the linear system for the coefficients $c_i$. Of course, it is also clear that each root $\xi_i(p)$ and coefficient $c_i$ have a companion conjugate root and conjugate coefficient, which implies that the imaginary parts appearing in the density cancel.
\end{remark}
\begin{remark}
In case that $\xi_i(p)$ has multiplicity greater than one for one or more values of $i$, the analysis proceeds in essentially the same way. For example, if  $\xi_1(p)=\xi_2(p)$, then the partial fraction decomposition of $\ltw$ becomes
$$
\ltw(s)=c_0+ \frac{c_1}{\bigl(s-\xi_1(p)\bigr)^2}+\sum_{i=2}^{N} \frac{c_i}{s-\xi_i(p)},
$$
the inverse of which is given by
$$
\dfw(x)={c}_1 x \mathrm{e}^{\xi_1(p) x}+\sum_{i=2}^{N} c_i \mathrm{e}^{\xi_i(p) x}.
$$
\end{remark}
\begin{remark}
For the nominal service time $B$ we have considered only mixtures of Erlang distributions, mainly because this class approximates well any continuous distribution on $[0, \infty)$ and because we can illustrate the techniques we use without complicating the analysis. However, we can extend this class by considering distributions with a rational Laplace transform. The analysis in \cite{vlasiou05} can be extended to such distributions, and the analysis in Cohen~\cite[Section II.5.10]{cohen-SSQ} is already given for such distributions, so the results given there can be implemented directly.
\end{remark}

\begin{remark}
The analysis we have presented so far can be directly extended to the case where $Y$ takes any finite number of negative values. In other words, let the distribution of $Y$ be given by $\p[Y=1]=p$, and for $i=1,\ldots,n$, $\p[Y=-u_i]=p_i$, where $u_i>0$ and $\sum_i p_i=1-p$. Then, for example, Equation~\eqref{eq:transform1} becomes
\begin{multline*}
\ltw(s)\Biggl[1-p\,\lta(-s)\sum_{n=1}^N \kappa_n\biggl(\frac{\mu}{\mu+s}\biggr)^n  \Biggr]\\=
p\,\p[W+B\leqslant A]-p\,\e[\mathrm{e}^{-s(B-A+W)};W+B\leqslant A]+\sum_{i=1}^n p_i\,\p[B\leqslant u_i W+A]+\\
+\sum_{i=1}^n p_i \Biggl[\sum_{n=1}^N \kappa_n \biggl(\frac{\mu}{\mu+s}\biggr)^n \lta(-s) \ltw(-u_i s)- \e[\mathrm{e}^{-s(B-A-u_i W)};B\leqslant u_i W+A]\Biggr].
\end{multline*}
Following the same steps as below \eqref{eq:transform1}, we can conclude that the waiting time density is again given by a mixture of exponentials of the form
$$
\dfw(x)=\sum_{i=1}^N \hat{c}_i \mathrm{e}^{\xi_i(p)x},
$$
where the new constants $\hat{c}_i$ (and the mass of the distribution at zero, given by $\hat{c}_0$) are to be determined as the unique solution to a linear system of equations. The only additional remark necessary when forming this linear system is to observe that both the probability $\p[B\leqslant u_i W+A]$ and the expectation $\e[\mathrm{e}^{-s(B-A-u_i W)};B\leqslant u_i W+A]$ can be expressed linearly in terms of the constants $\hat{c}_i$.
\end{remark}

\subsection*{The case $p=0$}
We have seen that the case where $Y_n=-1$ for all $n$, or in other words the case $p=0$, had to be excluded from the analysis. Equation~\eqref{eq:transform2} is still valid if we take the constants $\xi_i(0)$ to be defined as in Theorem~\ref{th:density}. However, one cannot apply Liouville's theorem to the resulting equation. The transform can be inverted directly. As it is shown in \cite{vlasiou05}, the terms $\phi^{(i)}(\mu)$ that remain to be determined follow by differentiating \eqref{eq:transform2} $N-1$ times and evaluating $\ltw[i](s)$ at $s=\mu$ for $i=0,\ldots,N-1$. The density in this case is a mixture of Erlang distributions with the same scale parameter $\mu$ for all exponential phases. As we can see, for $p=0$ the resulting density is intrinsically different from the one described in Theorem~\ref{th:density}.

\subsection*{The case $p=1$}
If $p=1$ and $\e[B]<\e[A]$, then we are analysing the steady-state waiting time distribution of a $\mathrm{G/PH/1}$ queue. Equation~\eqref{eq:transform2} now reduces to
\begin{multline}\label{eq:p=1 transform}
\ltw(s)\prod_{i=1}^{N}\bigl(s-\xi_i(1)\bigr)=\frac{\prod_{i=1}^{N}\bigl(s-\xi_i(1)\bigr)}{(\mu+s)^N-\lta(-s)\sum_{n=1}^N \kappa_n \mu^n(\mu+s)^{N-n}}\ \times\\
\times\Bigl[(\mu+s)^N\p[W+B\leqslant A]-(\mu+s)^N \e[\mathrm{e}^{-s(B-A+W)};W+B\leqslant A]\Bigl].
\end{multline}
Earlier we have already observed that the right-hand side of \eqref{eq:p=1 transform} is equal to an $N$-th degree polynomial $\sum_{i=0}^N q_i s^i$. Inspection of the right-hand side of \eqref{eq:p=1 transform} reveals that it has an $N$-fold zero in $s=-\mu$. Indeed, all zeros of the numerator of the quotient in the right-hand side cancel against zeros of the denominator, and the term
$$
\p[W+B\leqslant A]-\e[\mathrm{e}^{-s(B-A+W)};W+B\leqslant A]
$$
is finite for $s=-\mu$.
Hence,
\begin{equation}\label{eq:q_N}
\sum_{i=0}^N q_i s^i=q_N\, (\mu+s)^N.
\end{equation}
Combining \eqref{eq:p=1 transform} and \eqref{eq:q_N}, we conclude that
$$
\ltw(s)\prod_{i=1}^{N}\bigl(s-\xi_i(1)\bigr)=q_N\, (\mu+s)^N,
$$
and since $\ltw(0)=1$, the last equation gives us that
$$
q_N=\frac{\prod_{i=1}^{N}\bigl(-\xi_i(1)\bigr)}{\mu^N}.
$$
Thus, we have that
$$
\ltw(s)=\left(\frac{\mu+s}{\mu}\right)^N\prod_{i=1}^{N}\frac{\xi_i(1)}{\xi_i(1)-s},
$$
which is in agreement with Equation~II.5.190 in \cite[p.\ 324]{cohen-SSQ}.

\section{The M/D case}\label{s:M/D}
We have examined so far the case where the nominal interarrival time $A$ is generally distributed and the nominal service time $B$ follows a phase-type distribution. In other words, we have studied the case which is in a sense analogous to the ordinary $\mathrm{G/PH/1}$ queue. We now would like to study the reversed situation; namely, the case analogous to the $\mathrm{M/G/1}$ queue.

The $\mathrm{M/G/1}$ queue has been studied in much detail. However, the analogous alternating service model -- i.e.,
take $\p(Y=-1)=1$ in \eqref{eq:recursion}, so $p=0$ -- seems to be more complicated to analyse. As shown in \cite{vlasiou05a}, if $p=0$, the density of $W$ satisfies a generalised Wiener-Hopf equation, for which no solution is known in general. The presently available results for the distribution of $W$ with $p=0$ are developed in \cite{vlasiou05a}, where $B$ is assumed to belong to a class strictly bigger than the class of functions with rational Laplace transforms, but not completely general. Moreover, the method developed in \cite{vlasiou05a} breaks down when applied to \eqref{eq:equation} with $Y$ {\em not} identically equal to $-1$.

We shall refrain from trying to develop an alternative approach for the $\mathrm{M/G}$ case with a more general distribution for $B$ than the one treated in Section~\ref{s:G/PH}. Instead, we give a detailed analysis of the $\mathrm{M/D}$ case: $A$ is exponentially distributed and $B$ is deterministic. This case is neither contained in the $\mathrm{G/PH}$ case of the previous section nor has it been treated (for the special choice of $p=0$) in \cite{vlasiou05a}. Its analysis is of interest for various reasons. To start with, the model generalises the classical $\mathrm{M/D/1}$ queue; additionally, the analysis illustrates the difficulties that arise when studying \eqref{eq:equation} in case $A$ is exponentially distributed and $B$ is generally distributed; finally, the different
effects of Lindley's classical recursion and of the Lindley-type recursion discussed in \cite{vlasiou05a} are clearly exposed. As we shall see in the following, the analysis can be practically split into two parts, where each part follows the analysis of the corresponding model with $Y \equiv 1$, or $Y \equiv -1$.

\subsection{Deterministic nominal service times}
As before, consider Equation~\eqref{eq:equation}, and assume that $Y=1$ with probability $p$ and $Y=-1$ with probability $1-p$. Let $A$ be exponentially distributed with rate $\lambda$ and $B$ be equal to $b$, where $b>0$. Furthermore, we shall denote by $\pi_0$ the mass of the distribution of $W$ at zero; that is, $\pi_0=\p[W=0]$.

For this setting, we have from \eqref{eq:equation} that for $x\geqslant0$,
\begin{align}\label{eq:m/d distr}
\nonumber\Dfw(x)&=\p[\max\{0, b-A+YW\}\leqslant x]=\p[b-A+YW\leqslant x]\\
\nonumber        &=p\,\p[b-A+W\leqslant x]+(1-p)\,\p[b-A-W\leqslant x]\\
\nonumber        &=p\, \pi_0 \p[b-A\leqslant x]+p \int_0^\infty \p[b-A\leqslant x-y] \dfw(y) \,\mathrm{d}y+(1-p)\,\pi_0 \p[b-A\leqslant x]+\\
\nonumber        &\hspace{6cm}+(1-p) \int_0^\infty \p[b-A\leqslant x+y] \dfw(y) \,\mathrm{d}y\\
\nonumber        &=\pi_0 \p[A\geqslant b-x]+p \int_0^\infty \p[A\geqslant b-x+y] \dfw(y) \,\mathrm{d}y+\\
                 &\hspace{6cm}+(1-p) \int_0^\infty \p[A\geqslant b-x-y] \dfw(y)\,\mathrm{d}y.
\end{align}
So, for $0\leqslant x< b$ the above equation reduces to
\begin{multline}\label{eq:x<b distr}
\Dfw(x)=\pi_0\, \mathrm{e}^{-\lambda (b-x)}+p \int_0^\infty \mathrm{e}^{-\lambda (b-x+y)} \dfw(y) \,\mathrm{d}y+(1-p) \int_0^{b-x} \mathrm{e}^{-\lambda (b-x-y)} \dfw(y) \,\mathrm{d}y+\\
+(1-p) \int_{b-x}^\infty \dfw(y) \,\mathrm{d}y,
\end{multline}
and for $x\geqslant b$, Equation~\eqref{eq:m/d distr} reduces to
\begin{multline}\label{eq:x>b distr}
\Dfw(x)=\pi_0+p \int_0^{x-b}\dfw(y) \,\mathrm{d}y+p \int_{x-b}^\infty \mathrm{e}^{-\lambda (b-x+y)} \dfw(y) \,\mathrm{d}y+(1-p)(1-\pi_0),
\end{multline}
where we have utilised the normalisation equation
\begin{equation}\label{eq:normalisation}
\pi_0+\int_0^\infty\dfw(y) \,\mathrm{d}y=1.
\end{equation}

In the following, we shall derive the distribution on the interval $[0,b)$ and on the interval $[b,\infty)$ separately. At this point though, one should note that from Equation~\eqref{eq:equation} it is apparent that for $A$ exponentially distributed and $B=b$, the distribution of $W$ is continuous on $(0,\infty)$. Also, one can verify that Equation~\eqref{eq:x<b distr} for $x=b$ reduces to Equation~\eqref{eq:x>b distr} for $x=b$. The fact that $\Dfw$ is continuous on $(0,\infty)$ will be used extensively in the sequel. Notice also that from Equations~\eqref{eq:x<b distr} and \eqref{eq:x>b distr} we can immediately see that we can differentiate $\Dfw(x)$ for $x\in(0,b)$ and $x\in(b,\infty)$; see, for example, Titchmarsh~\cite[p.\ 59]{titchmarsh-tf}.\\

\noindent
\underline{\textbf{The distribution on} $\mathbf{[0,b)}.$}\\

In all subsequent equations it is assumed that $x\in(0,b)$. In order to derive the distribution of $W$ on $[0,b]$, we differentiate \eqref{eq:x<b distr} once to obtain
\begin{align*}
\dfw(x)&=\lambda \pi_0 \, \mathrm{e}^{-\lambda (b-x)}+\lambda p \int_0^\infty \mathrm{e}^{-\lambda (b-x+y)} \dfw(y) \,\mathrm{d}y+\lambda (1-p) \int_0^{b-x} \mathrm{e}^{-\lambda (b-x-y)} \dfw(y) \,\mathrm{d}y-\\
       &\quad -(1-p)  \mathrm{e}^{-\lambda (b-x)}  \mathrm{e}^{\lambda (b-x)}\dfw(b-x)+(1-p)\dfw(b-x).
\end{align*}
We rewrite this equation after noticing that the second line is equal to zero, while the sum of the integrals in the first line can be rewritten by using \eqref{eq:x<b distr}. Thus, we have that
\begin{align}\label{eq:x<b integral eq}
\nonumber \dfw(x)&=\lambda \pi_0 \, \mathrm{e}^{-\lambda (b-x)}+\lambda\left(\Dfw(x)-\pi_0 \mathrm{e}^{-\lambda (b-x)}-(1-p)\int_{b-x}^\infty \dfw(y) \,\mathrm{d}y\right)\\
                &=\lambda\Dfw(x)-\lambda(1-p)\int_{b-x}^\infty \dfw(y) \,\mathrm{d}y.
\end{align}
In order to obtain a linear differential equation, differentiate \eqref{eq:x<b integral eq} once more, which leads to
\begin{equation}\label{eq:x<b de}
\dfw['](x)=\lambda \dfw(x)-\lambda(1-p)\dfw(b-x).
\end{equation}

Equation \eqref{eq:x<b de} is a homogeneous linear differential equation, not of a standard form because of the argument $b-x$ that appears at the right-hand side. To solve it, we substitute $x$ for $b-x$ in \eqref{eq:x<b de} to obtain
\begin{equation}\label{eq:x<b de2}
\dfw['](b-x)=\lambda \dfw(b-x)-\lambda(1-p)\dfw(x).
\end{equation}
Then, we differentiate \eqref{eq:x<b de} once more to obtain
$$
\dfw[''](x)=\lambda \dfw['](x)+\lambda(1-p)\dfw['](b-x),
$$
and we eliminate the term $\dfw['](b-x)$ by using \eqref{eq:x<b de2}. Thus, we conclude that
\begin{equation}\label{eq:x<b main de}
\dfw[''](x)=\lambda^2 p\,(2-p) \dfw(x).
\end{equation}
For $p\neq 0$, the solution to this differential equation is given by
\begin{equation}\label{eq:x<b form of sol}
\dfw(x)=d_1 \mathrm{e}^{r_1 x}+d_2 \mathrm{e}^{r_2 x},
\end{equation}
where $r_1$ and $r_2$ are given by
\begin{equation}\label{eq:r12}
r_{1,2}=\pm \lambda \sqrt{p(2-p)},
\end{equation}
and the constants $d_1$ and $d_2$ will be determined by the initial conditions. Namely, the solution needs to satisfy \eqref{eq:x<b de} and the condition $\Dfw(0)=\pi_0$. Thus, for the first equation, substitute the general solution we have derived into \eqref{eq:x<b de}. For the second equation, first rewrite \eqref{eq:x<b integral eq} as follows:
$$
\dfw(x)=\lambda \Dfw(x)-\lambda (1-p) \left(1-\pi_0-\int_0^{b-x}\dfw(y)\,\mathrm{d}y\right),
$$
then substitute $\dfw(x)$ from \eqref{eq:x<b form of sol}, and finally evaluate the resulting equation for $x=0$. This system uniquely determines $d_1$ and $d_2$. Specifically, we have that
\begin{align*}
d_1&=\frac{\lambda^2 (1-p) \bigl(1-p\, (1-\pi_0)-2 \pi_0\bigr) r_1}{(\mathrm{e}^{b r_1}-1) \lambda^2 (2-p) (1-p)+\mathrm{e}^{b r_1} r_1 \bigl(r_1-\lambda (2-p)\bigr)},\\
d_2&=\frac{\mathrm{e}^{b r_1} \lambda (1-p\, (1-\pi_0)-2 \pi_0) r_1 \left(\lambda-r_1\right)}{(\mathrm{e}^{b r_1}-1) \lambda^2 (2-p) (1-p)+\mathrm{e}^{b r_1} r_1 \bigl(r_1-\lambda (2-p)\bigr)},
\end{align*}
where in the process we have assumed that $p\neq 1$. Up to this point we have that the waiting-time distribution on $[0,b]$ is given by
\begin{equation}\label{eq:distr x<b}
\Dfw(x)=\frac{d_1}{r_1}(\mathrm{e}^{r_1 x}-1)+\frac{d_2}{r_2}(\mathrm{e}^{r_2 x}-1)+\pi_0,
\end{equation}
where $d_1$ and $d_2$ are known up to the probability $\pi_0$. The cases for $p=0$ and $p=1$ follow directly from Equation~\eqref{eq:x<b main de} and will be handled separately in the sequel.\\

\noindent
\underline{\textbf{The distribution on} $\mathbf{[b,\infty)}.$}\\

As before, we obtain a differential equation by differentiating \eqref{eq:x>b distr} once, and substituting the resulting integrals by using \eqref{eq:x>b distr} once more. Thus, we obtain the equation
$$
\dfw(x)=\lambda\left(\Dfw(x)-\pi_0-(1-p)(1-\pi_0)-p\int_0^{x-b} \dfw(y)\,\mathrm{d}y\right),
$$
which can be reduced to
\begin{equation}\label{eq:x>b delay de}
\dfw(x)=\lambda\bigl(\Dfw(x)-1+p-p\,\Dfw(x-b)\bigr).
\end{equation}
Equation~\eqref{eq:x>b delay de} is a delay differential equation that can be solved recursively. Observe that for $x\in(b,2b)$, the term $\Dfw(x-b)$ has been derived in the previous step, so for $x\in(b,2b)$, Equation~\eqref{eq:x>b delay de} reduces to an ordinary linear differential equation from which we can easily derive the distribution of $W$ in the interval $(b,2b)$.

For simplicity, denote by $F_i (x)$ the distribution of $W$ when $x\in[ib,(i+1)b]$, and analogously denote by $f_i (x)$ the density of $W$, when $x\in(ib,(i+1)b)$. Then \eqref{eq:x>b delay de} states that
$$
f_i(x)=\lambda\bigl(F_i(x)-1+p-p\,F_{i-1}(x-b)\bigr),
$$
which leads to an expression for $F_i$ that is given in terms of an indefinite integral that is a function of $x$, that is,
\begin{equation}\label{eq:Distr formal x>b}
F_i(x)=\mathrm{e}^{\lambda x} \left[\int \lambda \bigl(-1+p-p\,F_{i-1}(x-b)\bigr) \mathrm{e}^{-\lambda x} \,\mathrm{d}x+\gamma_i\right],\quad i\geqslant 1.
\end{equation}
The constants $\gamma_i$ can be derived by exploiting the fact that the waiting-time distribution is continuous. In particular, every $\gamma_i$ is determined by the equation
\begin{equation}\label{8eq:asdf}
F_i(ib)=F_{i-1}(ib).
\end{equation}

Solving Equation~\eqref{eq:Distr formal x>b} recursively, we obtain that
\begin{multline}\label{distr M/D}
F_i(x)=1-p^i(1-\pi_0)-p^i\left(\frac{d_1}{r_1}+\frac{d_2}{r_2}\right)+\sum_{j=1}^2 \left(\frac{\lambda  p}{\lambda
-r_j}\right)^i \frac{d_j }{r_j} \mathrm{e}^{r_j (x-ib)}+\\
+x\sum_{j=0}^{i-1}(-\lambda p)^j \gamma_{i-j} \frac{(x-j b)^{j-1}}{j!} \mathrm{e}^{\lambda (x-jb)}.
\end{multline}
Observe that for $i=0$, if we define the empty sum at the right-hand side to be equal to zero, then the above expression is satisfied. Notice that, since we have made use of the distribution on $[0,b)$ as it is given by \eqref{eq:distr x<b}, Equation~\eqref{distr M/D} is not valid for $p=0$ or $p=1$. From Equation~\eqref{8eq:asdf} we now have that for every $i\geqslant 1$,
\begin{align}\label{eq:gammas}
\nonumber \gamma_i=&\,\mathrm{e}^{-\lambda i b}(1-p)p^{i-1}\left(\pi_0-1-\frac{d_1-d_2}{r_1}\right)-\sum_{j=1}^2 \frac{\mathrm{e}^{-\lambda i b}d_j }{r_j}\left(\frac{\lambda  p}{\lambda
-r_j}\right)^i\left(1-\frac{\mathrm{e}^{b r_j}\left(\lambda -r_j\right)}{\lambda  p}\right)+\\
&\quad+i \sum _{j=1}^{i-1} \frac{\mathrm{e}^{-\lambda  j b}(i-j)^{j-1}(-\lambda  p b)^j \left(\gamma _{i-1-j}-\gamma _{i-j}\right)}{j!}+\gamma_{i-1},
\end{align}
where we have assumed that $\gamma_0=0$, and that for $i=1$, the second sum is equal to zero. These expressions can be simplified further by observing that
$$
1-p^i(1-\pi_0)-p^i\left(\frac{d_1}{r_1}+\frac{d_2}{r_2}\right)=1-\frac{p^i}{2-p}.
$$

Recall that $d_1$ and $d_2$, and thus also all constants $\gamma_i$, are known in terms of $\pi_0$. The probability $\pi_0$ that still remains to be determined will be given by the normalisation equation \eqref{eq:normalisation}. Notice though, that since the waiting-time distribution is determined recursively for every interval $[ib,(i+1)b]$, Equation~\eqref{eq:normalisation} yields an infinite sum. The sum is well defined, since a unique density exists. The above findings are summarised in the following theorem.
\begin{theorem}\label{th:distr M/D}
Consider the recursion given by \eqref{eq:recursion}, and assume that $0<p<1$. Let $A$ be exponentially distributed with rate $\lambda$ and $B$ be equal to $b$, where $b>0$. Then for $x\in [ib,(i+1)b]$, $i=0,1,\ldots$, the limiting distribution of the waiting time is given by
\begin{multline*}
\Dfw(x)=1-\frac{p^i}{2-p}+\sum_{j=1}^2 \left(\frac{\lambda p}{\lambda -r_j}\right)^i \frac{d_j }{r_j} \mathrm{e}^{r_j (x-ib)}+x\sum_{j=0}^{i-1}(-\lambda p)^j \gamma_{i-j} \frac{(x-j b)^{j-1}}{j!} \mathrm{e}^{\lambda (x-jb)},
\end{multline*}
where the constants $\gamma_i$ are given by Equation~\eqref{eq:gammas} and the probability $\pi_0$ is given by the normalisation equation \eqref{eq:normalisation}.
\end{theorem}

\begin{figure}[ht]
\begin{center}
\includegraphics[width=0.8\textwidth]{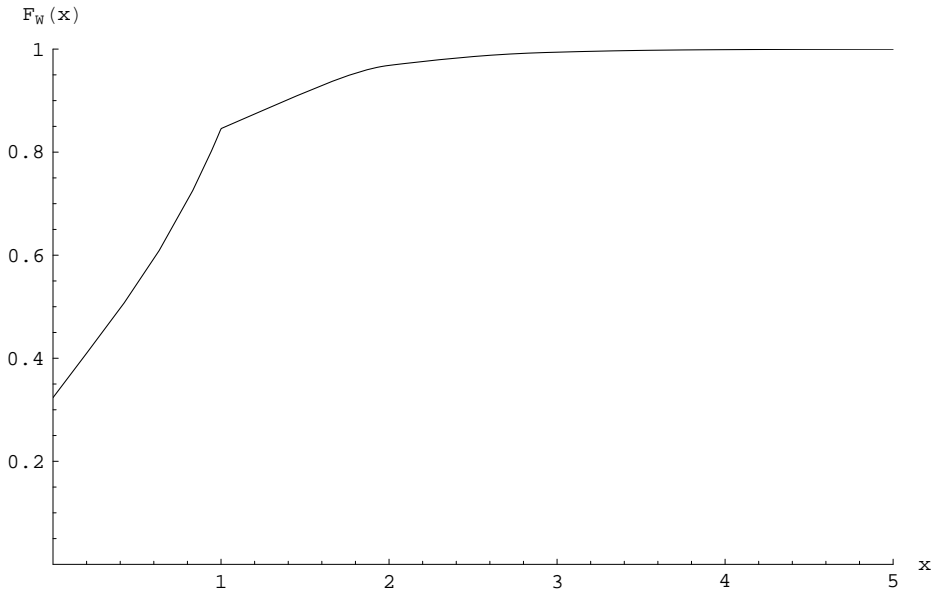}
\end{center}
\caption{The waiting time distribution for $b=1$, $\lambda=2$, and $p=1/3$.}
\label{fig:distr}
\end{figure}

One might expect though that Equation~\eqref{eq:normalisation} may not be suitable for numerically determining $\pi_0$. However, if the probability $p$ is not too close to one, or in other words, if the system does not almost behave like an $\mathrm{M/D/1}$ queue, then one can numerically approximate $\pi_0$ from the normalisation equation. As an example, in Figure~\ref{fig:distr} we display a typical plot of the waiting-time distribution. We have chosen $b=1$, $\lambda=2$, and $p=1/3$.

For $p$ close to one, we can see from the expressions for $d_1$ and $d_2$ that both the numerators and the denominators of these two constants approach zero. Furthermore, the denominators $\lambda-r_j$, $j=1,2$ that appear in the waiting-time distribution also approach zero, which makes Theorem~\ref{th:distr M/D} unsuitable for numerical computations for values of $p$ close to one. Moreover, we also see that very large values of the parameter $\lambda$ may also lead to numerical problems, since $\lambda$ is involved in the exponent of almost all exponential terms that appear in the waiting-time distribution.

As one can observe from Figure~\ref{fig:distr}, and show from Theorem~\ref{th:distr M/D}, $\Dfw$ is not differentiable for $x=b$. This is not surprising, as the waiting-time distribution is defined by two different equations; namely Equation~\eqref{eq:x<b distr} for $x<b$ and Equation~\eqref{eq:x>b distr} for $x\geqslant b$. Furthermore, from Equation~\eqref{eq:x<b integral eq} we have that
$$
\dfw(b^-)=\lambda \Dfw(b)-\lambda (1-p)(1-\pi_0),
$$
and from Equation~\eqref{eq:x>b delay de} we have that
$$
\dfw(b^+)=\lambda (\Dfw(b)-1+p-p \pi_0).
$$
That is, $\dfw(b^-)-\dfw(b^+)=\lambda \pi_0$.

\subsection*{ The case $p=0$}
Observe that if $p=0$ then the support of $W$ is the interval $[0,b]$. To determine the density of the waiting time, we insert $p=0$ into Equation~\eqref{eq:x<b main de}. Thus, we obtain that
$$
\dfw[''](x)=0,
$$
from which we immediately have that
$$
\dfw(x)=\nu_1x +\nu_2,
$$
for some constants $\nu_1$ and $\nu_2$ such that \eqref{eq:x<b de} is satisfied. The latter condition implies that for every $x \in(0,b)$ the following equation must hold:
$$
\nu_1=\lambda (\nu_1x +\nu_2)-\lambda (\nu_1 (b-x) +\nu_2).
$$
From this we conclude that $\nu_1$ is equal to zero, i.e.\ the waiting time has a mass at zero and is uniformly distributed on $(0,b)$. To determine the mass $\pi_0$ and the constant $\nu_2$ we evaluate \eqref{eq:x<b integral eq} at $x=0$ and we use the normalisation equation \eqref{eq:normalisation}, keeping in mind that $\dfw(x)=0$ for $x\in[b,\infty)$. These two equations yield that if $p=0$, then
\begin{equation}\label{eq:density p=0}
\dfw(x)=\frac{\lambda}{1+\lambda b},\quad 0<x<b, \qquad \mbox{and}\qquad \pi_0=\frac{1}{1+\lambda b}.
\end{equation}
Evidently, the density in this case is quite different from the density for $p\neq 0$, which is on $(0,b)$ a mixture of two exponentials; see \eqref{eq:x<b form of sol}.

Another way to see that $f_W(x) = \lambda \pi_0$, $0<x<b$, is as follows. Recall that for $p=0$ and $x\geqslant b$ we have that $\dfw(x)=0$. Equation~\eqref{eq:x<b integral eq} can now be written as
$$
f_W(x) = \lambda \pi_0 + \lambda \p[W \in (0,x)] - \lambda \p[W \in (b-x,b)].
$$
Replacing $x$ by $b-x$ shows that $\dfw(x) = \dfw(b-x)$, which implies that $\p[W \in (0,x)] = \p[W \in (b-x,b)]$ and finally that
$\dfw(x) = \lambda \pi_0$, $0<x<b$. It seems less straightforward to explain {\em probabilistically} that $W$, given that $W>0$, is uniformly distributed. With a view towards the recursion $W\stackrel{\m{D}}{=}\max\{0, b-A-W\}$, we believe that this property is related to the fact that, if $n$ Poisson arrivals occur in some interval, then they are distributed like the $n$ order statistics of the uniform distribution on that interval; see Ross~\cite[Section 2.3]{ross-SP}.

\subsection*{The case $p=1$}
For the $\mathrm{M/D/1}$ queue, Erlang~\cite{erlang09} derived the following expression for the waiting-time distribution:
\begin{equation*}
\p[W\leqslant x]=(1-\rho) \sum_{j=0}^i \frac{\bigl(-\lambda(x-jb)\bigr)^j}{j!} \, \mathrm{e}^{\lambda(x-jb)}, \quad ib\leqslant x < (i+1)b,
\end{equation*}
where $\rho$ is the traffic intensity. Recall that for the $\mathrm{M/D/1}$ queue we have that $\Dfw(0)=1-\rho$. We see that for $p=1$ Equation~\eqref{eq:x<b integral eq} indeed leads to the waiting-time distribution $(1-\rho)\, {\rm e}^{\lambda x}$, as it is given by Erlang's expression for the first interval $[0, b)$. For $x\geqslant b$, one needs to recursively solve Equation~\eqref{eq:Distr formal x>b} in order to obtain Erlang's expression. However, since the recursive solution we have obtained for our model makes use of $\Dfw(x)$ as it is given by \eqref{eq:distr x<b}, which is not valid for $p=1$, the waiting-time distribution we have obtained in Theorem~\ref{th:distr M/D} cannot be extended to the case for $p=1$.

The terms both in Erlang's expression for the waiting-time distribution of an $\mathrm{M/D/1}$ queue and in Theorem~\ref{th:distr M/D} alternate in sign and in general are much larger than their sum. Thus, the numerical evaluation of the sum may be hampered by roundoff errors due to the loss of significant digits, in particular under heavy traffic. For the $\mathrm{M/D/1}$ queue, however, a satisfactory solution has been given by Franx~\cite{franx01} in a way that only a finite sum of positive terms is involved; thus, this expression presents no numerical complications, not even for high traffic intensities. For our model, extending Franx's approach is a challenging problem as the representation of various quantities appearing in \cite{franx01} which are related to the queue length at service initiations is not straightforward.

As we see, the waiting-time distribution in Theorem~\ref{th:distr M/D} is quite similar to Erlang's expression, so we expect that eventually the solution will suffer from roundoff errors. However, a significant difference in the numerical computation between the $\mathrm{M/D/1}$ queue and the model described by Recursion~\eqref{eq:recursion} arises when computing $\pi_0$. For any single server queue we know a priori that $\p[W=0]=1-\rho$. In our model, $\pi_0$ has to be computed from the normalisation equation, where the numerical complications when calculating the waiting-time distribution become apparent. In particular, as $p$ tends to $1$, i.e.\ as the system behaves almost like an $\mathrm{M/D/1}$ queue, the computation of $\pi_0$ becomes more problematic.\\

As a final observation, we note that the effects of Lindley's classical recursion and of the Lindley-type recursion discussed in \cite{vlasiou05a} are quite apparent. The analysis for our model is in a sense separated into two parts: the derivation of the waiting-time distribution in $[0,b)$ and in $[b,\infty)$. In the first part, we see that Equation~\eqref{eq:x<b de} is quite similar to the differential equation appearing in \cite{vlasiou05b} for the derivation of the waiting-time distribution in case $p=0$ and $B$ follows a polynomial distribution. Moreover, one could use the same technique to derive a solution, but Equation~\eqref{eq:x<b de} is too simple to call for such means. In the second part, we see the effects of the $\mathrm{M/D/1}$ queue, as we eventually derive $\Dfw$ in a recursive manner. Furthermore, this model inherits all the numerical difficulties appearing in the classical solution for the $\mathrm{M/D/1}$ queue, plus the additional difficulties of computing $\pi_0$. For Lindley's recursion, $\pi_0$ is known beforehand, while for the Lindley-type recursion described in \cite{park03, vlasiou05a, vlasiou05, vlasiou05b, vlasiou04} $\pi_0$ is derived by the normalisation equation.

\section*{Acknowledgements}
The authors would like to thank Ivo Adan for his helpful remarks and his assistance in simulating this model.

\end{document}